\documentclass{article}
\usepackage{amssymb}
\usepackage{amsfonts}
\usepackage{amsmath}

\input{tcilatex}
\begin{document}

\begin{center}
\bigskip

\textbf{CLASSIFICATION RESULTS ON SURFACES IN THE ISOTROPIC 3-SPACE}

\bigskip

\textbf{Muhittin Evren AYDIN}

Department of Mathematics, Faculty of Science, Firat University, Elazig,
23119, Turkey meaydin@firat.edu.tr

\bigskip
\end{center}

\textbf{Abstract. }The isotropic 3-space $\mathbb{I}^{3}$ which is one of
the Cayley--Klein spaces is obtained from the Euclidean space by
substituting the usual Euclidean distance with the isotropic distance. In
the present paper, we give several classifications on the surfaces in $%
\mathbb{I}^{3}$ with the constant relative curvature (analogue of the
Gaussian curvature) and the constant isotropic mean curvature. In
particular, we classify the helicoidal surfaces in $\mathbb{I}^{3}$ with
constant curvature and analyze some special curves on these.

\bigskip

\textbf{Keywords:} Isotropic space; helicoidal surface; isotropic mean
curvature; relative curvature.

\textbf{Math. Subject Classification 2010: }$53A35$, $53A40$, $53B25$.

\section{Preliminaries}

Differential geometry of isotropic spaces have been introduced by K.
Strubecker \cite{41}, H. Sachs \cite{34}-\cite{36}, D. Palman \cite{29} and
others. Especially the reader can find a well bibliography for isotropic
planes and isotropic 3-spaces in \cite{34,35}.

The isotropic 3-space $\mathbb{I}^{3}$ is a Cayley--Klein space defined from
a 3-dimensional projective space $P\left( 
\mathbb{R}
^{3}\right) $ with the absolute figure which is an ordered triple $\left(
\omega ,f_{1},f_{2}\right) $, where $\omega $ is a plane in $P\left( 
\mathbb{R}
^{3}\right) $ and $f_{1},f_{2}$ are two complex-conjugate straight lines in $%
\omega $ (see \cite{39}). The homogeneous coordinates in $P\left( 
\mathbb{R}
^{3}\right) $ are introduced in such a way that the \textit{absolute plane} $%
\omega $ is given by $X_{0}=0$ and the \textit{absolute lines} $f_{1},f_{2}$
by $X_{0}=X_{1}+iX_{2}=0,$ $X_{0}=X_{1}-iX_{2}=0.$ The intersection point $%
F(0:0:0:1)$ of these two lines is called the \textit{absolute point}. The
group of motions of $\mathbb{I}^{3}$ is a six-parameter group given in the
affine coordinates $x_{1}=\frac{X_{1}}{X_{0}},$ $x_{2}=\frac{X_{2}}{X_{0}},$ 
$x_{3}=\frac{X_{3}}{X_{0}}$ by

\begin{equation}
\left( x_{1},x_{2},x_{3}\right) \longmapsto \left( x_{1}^{\prime
},x_{2}^{\prime },x_{3}^{\prime }\right) :\left\{ 
\begin{array}{l}
x_{1}^{\prime }=a+x_{1}\cos \phi -x_{2}\sin \phi , \\ 
x_{2}^{\prime }=b+x_{1}\sin \phi +x_{2}\cos \phi , \\ 
x_{3}^{\prime }=c+dx_{1}+ex_{2}+x_{3},%
\end{array}%
\right.  \tag{1.1}
\end{equation}%
where $a,b,c,d,e,\phi \in 
\mathbb{R}
.$

Such affine transformations are called \textit{isotropic congruence
transformations }or \textit{i-motions. }It easily seen from $\left(
1.1\right) $ that i-motions are indeed composed by an Euclidean motion in
the $x_{1}x_{2}-$plane (i.e. translation and rotation) and an affine shear
transformation in $x_{3}-$direction.

Consider the points $\mathbf{x}=\left( x_{1},x_{2},x_{3}\right) $ and $%
\mathbf{y}=\left( y_{1},y_{2},y_{3}\right) .$ The projection in $x_{3}-$%
direction onto $\mathbb{R}^{2},$ $\left( x_{1},x_{2},x_{3}\right)
\longmapsto \left( x_{1},x_{2},0\right) ,$ is called the\textit{\ top view}.
The \textit{isotropic distance}, so-called \textit{i-distance }of two points 
$\mathbf{x}$ and $\mathbf{y}$ is defined as the Euclidean distance of their
top views, i.e.,%
\begin{equation}
\left\Vert \mathbf{x}-\mathbf{y}\right\Vert _{i}=\sqrt{\sum_{j=1}^{2}\left(
y_{j}-x_{j}\right) ^{2}}.  \tag{1.2}
\end{equation}%
The i-metric is degenerate along the lines in $x_{3}-$direction, and such
lines are called \textit{isotropic} lines.

\textit{Planes, circles and spheres}. There are two types of planes in $%
\mathbb{I}^{3}$ (\cite{31}-\cite{33}).

(1) \textit{Non-isotropic planes} are planes non-parallel to the $x_{3}-$%
direction. In these planes we basically have an Euclidean metric: This is
not the one we are used to, since we have to make the usual Euclidean
measurements in the top view. An \textit{i-circle} (of \textit{elliptic type}%
) in a non-isotropic plane $P$ is an ellipse, whose top view is an Euclidean
circle. Such an i-circle with center $\mathbf{m}\in P$ and radius $r$ is the
set of all points\textbf{\ }$\mathbf{x}\in P$ with $\left\Vert \mathbf{x}-%
\mathbf{m}\right\Vert _{i}=r.$

(2) \textit{Isotropic planes }are planes parallel to the $x_{3}-$axis.
There, $\mathbb{I}^{3}$ induces an isotropic metric. An \textit{i-circle }%
(of \textit{parabolic type}) is a parabola with $x_{3}-$parallel axis and
thus it lies in an isotropic plane

An i-circle of parabolic type is not the iso-distance set of a fixed point,
but it may be seen as a curve with constant isotropic curvature: A curve $c$
in an isotropic plane P (without loss of generality we set $P:x_{2}=0$)
which does not possess isotropic tangents can be written as graph $%
x_{3}=f(x_{1})$. Then, the \textit{i-curvature} of $c$ at $x_{1}=\alpha _{0}$
is given by the second derivative $\kappa _{i}\left( \alpha _{0}\right)
=f^{\prime \prime }\left( \alpha _{0}\right) $. For an i-circle of parabolic
type $f$ is quadratic and thus $\kappa _{i}$ is constant.

There are also two types of \textit{isotropic spheres}. An \textit{i-sphere
of the cylindrical type} is the set of all points $\mathbf{x}\in \mathbb{I}%
^{3}$ with $\left\Vert \mathbf{x}-\mathbf{m}\right\Vert _{i}=r$. Speaking in
an Euclidean way, such a sphere is a right circular cylinder with $x_{3}$%
\textit{-}parallel rulings; its top view is the Euclidean circle with center 
$\mathbf{m}$ and radius $r$. The more interesting and important type of
spheres are the \textit{i-spheres of parabolic type,}%
\begin{equation*}
x_{3}=\frac{A}{2}\left( x_{1}^{2}+x_{2}^{2}\right) +Bx_{1}+Cx_{2}+D,\text{ \
\ }A\neq 0.
\end{equation*}

From an Euclidean perspective, they are paraboloids of revolution with $z-$%
parallel axis. The intersections of these i-spheres with planes $P$ are
i-circles. If $P$ is non-isotropic, then the intersection is an i-circle of
elliptic type. If $P$ is isotropic, the intersection curve is an i-circle of
parabolic type.

\textit{Curvature theory of surfaces. }A surface $M^{2}$ immersed in $%
\mathbb{I}^{3}$ is called \textit{admissible} if it has no isotropic tangent
planes. We restrict our framework to admissible regular surfaces. For such a
surface $M^{2}$, the coefficients $g_{11}$, $g_{12}$, $g_{22}$ of its first
fundamental form are calculated with respect to the induced metric.

The normal field of $M^{2}$ is always the isotropic vector $\left(
0,0,1\right) $ since it is perpendicular to all tangent vectors to $M^{2}.$
The coefficients $h_{11}$, $h_{12}$, $h_{22}$ of the second fundamental form
of $M^{2}$ are calculated with respect to the normal field of $M^{2}$ (for
details, see \cite{35}, p. 155).

The \textit{relative curvature }(so called \textit{isotropic Gaussian
curvature}) and \textit{isotropic mean curvature} are defined by%
\begin{equation*}
K=\frac{\det \left( h_{ij}\right) }{\det \left( g_{ij}\right) },\text{ \ \ }%
H=\frac{g_{11}h_{22}-2g_{12}h_{12}+g_{22}h_{11}}{2\det \left( g_{ij}\right) }%
.
\end{equation*}

The surface $M^{2}$ is said to be \textit{isotropic flat} (resp. \textit{%
isotropic minimal}) if $K$ (resp.$H)$ vanishes.

The surfaces in the isotropic spaces have been studied by I. Kamenarovic (%
\cite{21,22}), B. Pavkovic (\cite{30}) and B. Divjak (\cite{14,40}) as well
as has important applications in several applied sciences, e.g., computer
science, Image Processing, architectural design and microeconomics (\cite%
{9,12,26}, \cite{31}-\cite{33}).

Most recently, Z. M. Sipus (\cite{39})\ classified the translation surfaces
of constant curvature generated by two planar curves in $\mathbb{I}^{3}$.
And then a classification for the ones generated by a space curve and a
planar curve with constant curvature was derived in \cite{2}.

In \cite{3}, the authors established a method to calculate the second
fundamental form of the surfaces in the isotropic 4-space $\mathbb{I}^{4}$
and classified some surfaces in $\mathbb{I}^{4}$ with vanishing curvatures.

In this paper, the helicoidal surfaces in $\mathbb{I}^{3}$ with constant
isotropic mean and constant relative curvature are classified. Further some
special curves on such surfaces are characterized.

\section{Classification results on surfaces in $\mathbb{I}^{3}$}

This section is devoted to recall the classification results on surfaces in $%
\mathbb{I}^{3}$ into seperate subsections, such as the translation surfaces,
the homothetical surfaces (so-called factorable surfaces), Aminov surfaces,
the spherical product surfaces.

\subsection{Translation surfaces in $\mathbb{I}^{3}$}

The present author introduced the translation surfaces generated by a space
curve and a planar curve as follows (for details, see \cite{2})%
\begin{equation}
\mathbf{r}\left( u_{1},u_{2}\right) \longmapsto \left( f_{1}\left(
u_{1}\right) ,f_{2}\left( u_{1}\right) +g_{2}\left( u_{2}\right)
,f_{3}\left( u_{1}\right) +g_{3}\left( u_{2}\right) \right) ,  \tag{2.1}
\end{equation}%
and classified the ones with constant curvature by the following theorems:

\bigskip

\textbf{Theorem 2.1.} \cite{2} \textit{Let }$M^{2}$\textit{\ be a
translation surface given by }$\left( 2.1\right) $\textit{\ in }$\mathbb{I}%
^{3}$\textit{\ with constant relative curvature }$K_{0}$\textit{. Then it is
either a generalized cylinder, i.e. }$K_{0}=0,$ \textit{or parametrized by
one of the following}

(i) $\mathbf{r}\left( u_{1},u_{2}\right) =\left(
f_{1},b_{1}f_{1}+g_{2},a_{1}\left( f_{1}\right) ^{2}+\frac{K_{0}}{a_{1}}%
\left( g_{2}\right) ^{2}+b_{2}f_{1}+b_{3}g_{2}\right) ;$

(ii) $\mathbf{r}\left( u_{1},u_{2}\right) =\left( f_{1},a_{2}\left(
f_{1}\right) ^{2}+b_{4}f_{1}+g_{2},b_{5}\left( f_{1}\right) ^{2}+\frac{1}{%
K_{0}a_{3}}\left( -2K_{0}g_{2}\right) ^{3/2}+b_{6}f_{1}+a_{4}g_{2}\right) ,$%
\newline
\textit{where }$a_{i}$\textit{\ are nonzero constants and }$b_{j}$\textit{\
some constants for }$1\leq i\leq 4$\textit{\ and }$1\leq j\leq 6.$

\bigskip

\textbf{Theorem 2.2.} \cite{2} \textit{Let }$M^{2}$\textit{\ be a
translation surface given by }$\left( 2.1\right) $\textit{\ in }$\mathbb{I}%
^{3}$ \textit{with constant isotropic mean curvature }$H_{0}.$\textit{\ Then
it is determined by one of the following expressions}

(i) $\mathbf{r}\left( u_{1},u_{2}\right) =\left(
f_{1},f_{2}+g_{2},H_{0}f_{1}^{2}+b_{1}f_{2}+b_{2}f_{1}+b_{3}g_{2}\right) ,$

(ii) $\mathbf{r}\left( u_{1},u_{2}\right) =\left(
f_{1},b_{4}f_{1}+g_{2},\left( H_{0}-a_{1}\right) \left( f_{1}\right)
^{2}+a_{2}\left( g_{2}\right) ^{2}+b_{5}f_{1}+b_{6}g_{2}\right) ;$

(iii) $\mathbf{r}\left( u_{1},u_{2}\right) =\left( f_{1},-\frac{1}{a_{3}}\ln
\left\vert \cos \left( a_{3}f_{1}\right) \right\vert
+g_{2},H_{0}f_{1}^{2}+b_{7}f_{2}+\frac{1}{a_{3}^{2}}\exp \left(
a_{3}g_{2}\right) +b_{8}f_{1}+b_{9}g_{2}\right) ,$\newline
\textit{where }$a_{i}$\textit{\ are nonzero constant and }$b_{i}$\textit{\
some constants }$1\leq i\leq 3$ and $1\leq j\leq 9.$

\bigskip

\textbf{Remark 2.3. }\textit{Isotropic minimal translation surfaces can also
be classified by Theorem 2.2 as taking }$H_{0}=0$\textit{\ in the statements
(i)-(iii) of the theorem.}

\bigskip

The author and A.O. Ogrenmis (\cite{4}) introduced the translation
hypersurfaces in such a that way a \textit{translation hypersurface} $\left(
M^{n},F\right) $ in the isotropic $\left( n+1\right) -$space $\mathbb{I}%
^{n+1},$ $n\geq 2,$ is parametrized by%
\begin{equation*}
X:\mathbb{R}^{n}\longrightarrow \mathbb{I}^{n+1},\text{ }\mathbf{x}%
\longmapsto \left( \mathbf{x},F\left( \mathbf{x}\right) \right) ,\text{ }%
F\left( \mathbf{x}\right) :=\sum_{j=1}^{n}f_{j}\left( x_{j}\right) ,\text{ }%
\mathbf{x}\in \mathbb{R}^{n},
\end{equation*}%
where $f_{j}$ is a smooth function of one variable for all $j\in \left\{
1,...,n\right\} .$ For more details of $\mathbb{I}^{n+1},$ see \cite{9,36,40}%
.

Some classifications were obtained for such hypersurfaces in $\mathbb{I}%
^{n+1}$ by the following results:

\bigskip

\textbf{Theorem 2.4. }\cite{4}\textbf{\ }\textit{Let} $\left( M^{n},F\right) 
$ \textit{be a translation hypersurface in} $\mathbb{I}^{n+1}$\textit{\ with
nonzero constant relative curvature }$K_{0}$\textit{. Then it has of the form%
}%
\begin{equation*}
X\left( \mathbf{x}\right) =\left( \mathbf{x},\sum_{j=1}^{n}\alpha
_{j}x_{j}^{2}+\beta _{j}x_{j}+\varepsilon \right) ,
\end{equation*}%
\textit{where }$\mathbf{x}\in \mathbb{R}^{n}$,\textit{\ }$\alpha _{j}$ are
nonzero constants and $\beta _{j},\varepsilon $\textit{\ some constants for
all }$j\in \left\{ 1,...,n\right\} .$\textit{\ }

\textit{In particular, if }$\left( M^{n},F\right) $ \textit{is isotropic
flat in} $\mathbb{I}^{n+1},$ \textit{then it\ is congruent to a cylinder
from Euclidean perspective.}

\bigskip

\textbf{Theorem 2.5. }\cite{4}\textbf{\ }\textit{Let} $\left( M^{n},F\right) 
$ \textit{be a translation hypersurface in} $\mathbb{I}^{n+1}$\textit{\ with
constant isotropic mean curvature }$H_{0}$\textit{. Then it has of the form}%
\begin{equation}
X\left( \mathbf{x}\right) =\left( \mathbf{x},\sum_{j=1}^{n}\alpha
_{j}x_{j}^{2}+\beta _{j}x_{j}+\varepsilon \right) ,  \tag{4.1}
\end{equation}%
\textit{where }$\mathbf{x}\in \mathbb{R}^{n}$\textit{\ and }$\alpha
_{j},\beta _{j},\varepsilon $\textit{\ are some constants for all }$j\in
\left\{ 1,...,n\right\} $\textit{\ such that }$\sum_{j=n}^{n}\alpha _{j}=%
\frac{n}{2}H_{0}.$\textit{\ }

\bigskip

\textbf{Remark 2.6. }\textit{Isotropic minimal translation hypersurfaces in} 
$\mathbb{I}^{n+1}$ \textit{are also classified by Theorem 2.5 as taking} $%
H_{0}=0.$

\subsection{Homothetical surfaces in $\mathbb{I}^{n+1}$}

The authors in \cite{4} defined the homothetical hypersurfaces in $\mathbb{I}%
^{n+1}$ as follows: A hypersurface $M^{n}$ of $\mathbb{I}^{n+1}$ is called a 
\textit{homothetical hypersurface }$\left( M^{n},H\right) $ if it is the
graph of a function of the form: 
\begin{equation*}
H\left( x_{1},...,x_{n}\right) :=h_{1}\left( x_{1}\right) \cdot ...\cdot
h_{n}\left( x_{n}\right) ,
\end{equation*}%
where $h_{1},...,h_{n}$ are smooth functions of one real variable.

Next results classify the homothetical hypersurfaces in $\mathbb{I}^{n+1}$
with constant isotropic mean and relative curvature.

\bigskip

\textbf{Theorem 2.7. }\cite{4} \textit{Let }$\left( M^{n},H\right) $ \textit{%
be a homothetical hypersurface in }$\mathbb{I}^{n+1}$ \textit{with constant
isotropic mean curvature }$H_{0}$. \textit{Then it is isotropic minimal,
i.e. }$H_{0}=0$\textit{\ and has the following form}

\begin{equation}
X\left( \mathbf{x}\right) =\left( \mathbf{x},\dprod\limits_{j=1}^{n}\left\{
\gamma _{j}x_{j}+\varepsilon _{j}\right\} \right) ,  \tag{3.1}
\end{equation}%
\textit{where }$\mathbf{x\in }\mathbb{R}^{n},$\textit{\ }$\gamma _{j}$%
\textit{, }$\varepsilon _{j}$\textit{\ some constants.}

\bigskip

\textbf{Theorem\textbf{\ }2.8.} \cite{4}\textbf{\ }\textit{Let }$\left(
M^{n},H\right) $\textit{\ be an isotropic flat homothetical hypersurface in }%
$\mathbb{I}^{n+1}$\textit{. Then it has one of the following forms:}

(ii) $X\left( \mathbf{x}\right) =\left( \mathbf{x},\gamma \exp \left( \alpha
_{1}x_{1}+\alpha _{2}x_{2}\right) \tprod\limits_{j=3}^{n}h_{j}\left(
x_{j}\right) \right) \mathit{\ }$\textit{for nonzero constants }$\gamma
,\alpha _{1},\alpha _{2};$

(ii)\textit{\ }$X\left( \mathbf{x}\right) =\left( \mathbf{x},\gamma
\tprod\limits_{j=1}^{n}\left( x_{j}+\beta _{j}\right) ^{\alpha _{j}}\right)
, $ \textit{where }$\mathbf{x}\in \mathbb{R}^{n},$\textit{\ }$\beta
_{1},...,\beta _{n}$\textit{\ are some constants and }$\gamma ,\alpha
_{1},...,\alpha _{n}$\textit{\ nonzero constants such that }$%
\sum_{i=1}^{n}\alpha _{i}=1.$

\bigskip

In the particular three dimensional case, the same authors generalized the
Theorem 2.8 to the homothetical surfaces with any constant relative curvature

\bigskip

\textbf{Theorem 2.9.} \cite{5}\textbf{\ }\textit{Let }$M^{2}$\textit{\ be a
homothetical surface\ in }$\mathbb{I}^{3}$ \textit{with constant relative
curvature }$K_{0}$.

\textit{(A) If }$K_{0}=0,$\textit{\ then we have }

\textit{(A.1)\ }$H\left( x,y\right) =c_{1}h_{1}\left( x\right) $\textit{\ or 
}$H\left( x,y\right) =c_{2}h_{2}\left( y\right) $\textit{\ for nonzero
constants }$c_{1},c_{2},$

\textit{(A.2)\ }$H\left( x,y\right) =c_{1}\exp \left( c_{2}x+c_{3}y\right) ,$%
\textit{\ where }$c_{1},c_{2},c_{3}$\textit{\ are nonzero constants,}

\textit{(A.3)\ }$H\left( x,y\right) =\left( c_{1}x+d_{1}\right) ^{\frac{1}{%
c_{2}}}\left( c_{3}y+d_{2}\right) ^{\frac{1}{c_{4}}},$\textit{\ where }$%
c_{1},c_{2},c_{3},c_{4}$\textit{\ are nonzero constants and }$d_{1},d_{2}$%
\textit{\ some constants, }$c_{2}\neq 1\neq c_{4}$\textit{, }

\textit{(B) If }$K_{0}\neq 0,$ \textit{then it is negative (i.e. }$K_{0}<0$%
\textit{) and }$h_{1}$, $h_{2}$\textit{\ are linear functions.}

\subsection{Spherical product surfaces and Aminov surfaces in $\mathbb{I}%
^{4} $}

The present author and I. Mihai (see \cite{3}) established a method to
calculate the second fundamental form of the surfaces in the isotropic
4-space $\mathbb{I}^{4}.$ Then ones classified the Aminov surfaces, given by%
\begin{equation*}
\mathbf{r}:I\times \left[ 0,2\pi \right) \longrightarrow \mathbb{I}^{4},%
\text{ }\left( u,v\right) \longmapsto \mathbf{r}\left( u,v\right) :=\left(
u,v,r\left( u\right) \cos v,r\left( u\right) \sin v\right) ,
\end{equation*}%
with vanishing curvature as follows:

\bigskip

\textbf{Theorem 2.10. }\cite{3}\textbf{\ }\textit{The isotropic flat Aminov
surfaces in }$\mathbb{I}^{4}$\textit{\ are only generalized cylinders over
circular helices from Euclidean perspective.}

\bigskip

\textbf{Theorem 2.11. }\cite{3}\textbf{\ }\textit{There does not exist an
isotropic minimal Aminov surface in }$\mathbb{I}^{4}.$

\bigskip

Furthermore, same authors derived the following classification results for
the spherical product surface $\left( M^{2}\mathit{,}c_{1}\otimes
c_{2}\right) $ of two curves $c_{1}$ and $c_{2}$ in $\mathbb{I}^{4}$ which
is defined by%
\begin{equation*}
\mathbf{r:}=c_{1}\otimes c_{2}:\mathbb{R}^{2}\longrightarrow \mathbb{I}^{4},%
\text{ }\left( u,v\right) \longmapsto \left( u,f_{1}\left( u\right)
,f_{2}\left( u\right) v,f_{2}\left( u\right) g\left( v\right) \right) ,
\end{equation*}%
where the curves $c_{1}\left( u\right) =\left( u,f_{1}\left( u\right)
,f_{2}\left( u\right) \right) $ and $c_{2}\left( v\right) =\left( v,g\left(
v\right) \right) $ are called the \textit{generating curves} of the surface.

\bigskip

\textbf{Theorem 2.12. }\cite{3}\textbf{\ }\textit{Let }$\left( M^{2}\mathit{,%
}c_{1}\otimes c_{2}\right) $\textit{\ be a isotropic flat spherical product
surface in} $\mathbb{I}^{4}.$ \textit{Then either it is a non-isotropic
plane or one of the following satisfies}

\textit{(i) }$c_{1}$ \textit{is a planar curve in }$\mathbb{I}^{3}$ \textit{%
lying in the isotropic plane }$y=const.;$

\textit{(ii) }$c_{1}$ \textit{is a line in }$\mathbb{I}^{3};$

\textit{(iii)} $c_{1}$ is\textit{\ a curve in }$\mathbb{I}^{3}$ \textit{of
the form }%
\begin{equation*}
c_{1}\left( u\right) =\left( u,f_{1}\left( u\right) ,\lambda \int \sqrt{%
1+\left( f_{1}^{\prime }\right) ^{2}}du+\xi \right) ,\text{ }\lambda ,\xi
\in 
\mathbb{R}
,\text{ }\lambda \neq 0;
\end{equation*}

(iv) $c_{2}$\textit{\ is a line in }$\mathbb{I}^{2}.$

\bigskip

\textbf{Theorem 2.13. }\cite{3}\textbf{\ }\textit{There does not exist an is
isotropic minimal spherical product surface in} $\mathbb{I}^{4}$ \textit{%
excepting totally geodesic ones.}

\section{Helicoidal surfaces in $\mathbb{I}^{3}$}

Rotation surfaces in the Euclidean 3-space $\mathbb{R}^{3}$ with constant
mean curvature have been known for a long time \cite{13,23}. A natural
generalization of rotation surfaces in $\mathbb{R}^{3}$ are the \textit{%
helicoidal surfaces }that can be defined as the orbit of a plane curve under
a screw motion in $\mathbb{R}^{3}$.

Such surfaces in $\mathbb{R}^{3}$ with constant mean and Gaussian curvature
have been classified by M. Do Carmo and M. Dajczer in \cite{15}. These
classifications were extended to the ones with prescribed mean and Gaussian
curvatures by C. Baikoussis and T. Koufogiorgos \cite{7}.

The helicoidal surfaces also have been studied by many authors in the
Minkowskian 3-space $\mathbb{R}_{1}^{3}$, the pseudo-Galilean space $\mathbb{%
G}_{3}^{1}$ and several homogeneous spaces as focusing on curvature
properties (see \cite{1,8,14,19,20,25,27,28}).

Morever, there exist many works related with the helicoidal surfaces
satisfying an equation in terms of its position vector and Laplace operator
in $\mathbb{R}^{3}$ and $\mathbb{R}_{1}^{3}$. For example see \cite%
{6,11,24,37,38}.

Now we adapt the above notion to isotropic spaces. Considering the i-motions
given by $\left( 1.1\right) ,$ the Euclidean rotations in the isotropic
space $\mathbb{I}^{3}$ is given by the normal form (in affine coordinates)%
\begin{equation*}
\left\{ 
\begin{array}{l}
x_{1}^{\prime }=x_{1}\cos \phi -x_{2}\sin \phi , \\ 
x_{2}^{\prime }=x_{1}\sin \phi +x_{2}\cos \phi , \\ 
x_{3}^{\prime }=x_{3},%
\end{array}%
\right.
\end{equation*}%
where $\phi \in 
\mathbb{R}
.$

Now let $c$ be a curve lying in the isotropic $x_{1}x_{3}-$plane given by $%
c\left( u\right) =\left( f\left( u\right) ,0,g\left( u\right) \right) $
where $f,g\in C^{2}$ and $f\neq 0\neq \frac{df}{du}$. By rotating the curve $%
c$ around $z-$axis and simultaneously followed by a translation, we obtain
that the \textit{helicoidal surface of first type }in $\mathbb{I}^{3}$ with
the profil curve $c$ and pitch $h$ is of the form

\begin{equation}
\mathbf{r}\left( u,v\right) =\left( f\left( u\right) \cos v,f\left( u\right)
\sin v,g\left( u\right) +hv\right) ,\text{ }h\in 
\mathbb{R}
.  \tag{2.1}
\end{equation}%
Similarly when the profile curve $c$ lies in the isotropic $x_{2}x_{3}-$%
plane, then the\textit{\ helicoidal surface} \textit{of second type }in $%
\mathbb{I}^{3}$ with pitch $h$ is given by

\begin{equation}
\mathbf{r}\left( u,v\right) =\left( -f\left( u\right) \sin v,f\left(
u\right) \cos v,g\left( u\right) +hv\right) ,\text{ }h\in 
\mathbb{R}
.  \tag{2.2}
\end{equation}%
In case $h=0,$ these reduce to the \textit{surfaces of revolution} in $%
\mathbb{I}^{3}.$ Also when $g$ is a constant, then $\mathbf{r}$ is a \textit{%
helicoid} from Euclidean perspective.

\bigskip

\textbf{Remark 2.1. }The coordinate functions $f$ and $g$ of the profile
curve $c$ are arbitrary functions of class $C^{2}$ and so one can take $%
f\left( u\right) =u.$

\bigskip

\textbf{Remark 2.2. }Since both type of the helicoidal surfaces are locally
isometric, we only will focus on the ones of first type.

\bigskip

Let $M^{2}$ be a helicoidal surface of first type in $\mathbb{I}^{3}.$ Then
the matrix of the first fundamental form $\mathfrak{g}$ of $M^{2}$ is%
\begin{equation*}
\left( \mathfrak{g}_{ij}\right) =%
\begin{bmatrix}
1 & 0 \\ 
0 & u^{2}%
\end{bmatrix}%
\text{ and }\left( \mathfrak{g}^{ij}\right) =%
\begin{bmatrix}
1 & 0 \\ 
0 & \frac{1}{u^{2}}%
\end{bmatrix}%
,
\end{equation*}%
where the prime denotes the derivative with respect to $u$ and $\left( 
\mathfrak{g}^{ij}\right) =\left( \mathfrak{g}_{ij}\right) ^{-1}.$ Thus the
Laplacian of $M^{2}$ with respect to $\mathfrak{g}$ is%
\begin{equation*}
\bigtriangleup =\frac{1}{\sqrt{\det \left( \mathfrak{g}_{ij}\right) }}%
\sum_{i.j=1}^{2}\frac{\partial }{\partial u_{i}}\left( \sqrt{\det \left( 
\mathfrak{g}_{ij}\right) }\mathfrak{g}^{ij}\frac{\partial }{\partial u_{j}}%
\right)
\end{equation*}%
and by taking $u_{1}=u$ and $u_{2}=v,$\ we get%
\begin{equation*}
\bigtriangleup =\frac{1}{u}\frac{\partial }{\partial u}+\frac{\partial ^{2}}{%
\partial u^{2}}+\frac{1}{u^{2}}\frac{\partial ^{2}}{\partial v^{2}}.
\end{equation*}%
Putting $r_{1}\left( u,v\right) =u\cos v$, $r_{2}\left( u,v\right) =u\sin v$
and $r_{3}\left( u,v\right) =g\left( u\right) +hv,$ one can easily seen that 
$\bigtriangleup r_{i}=0,$ $i=1,2$ and 
\begin{equation*}
\bigtriangleup r_{3}=\frac{1}{u}g^{\prime }+g^{\prime \prime }.
\end{equation*}%
Assuming $\bigtriangleup r_{3}=\lambda r_{3},$ $\lambda \in 
\mathbb{R}
,$ we can obtain that $\lambda $ must be zero and 
\begin{equation}
\frac{1}{u}g^{\prime }+g^{\prime \prime }=0.  \tag{2.3}
\end{equation}%
After solving $\left( 2.3\right) ,$ we derive $g\left( u\right) =\alpha \ln
\left\vert u\right\vert +\beta $ for $\alpha \in 
\mathbb{R}
\backslash \left\{ 0\right\} ,$ $\beta \in 
\mathbb{R}
.$

Thus we have the following result

\bigskip

\textbf{Proposition 2.2.} \textit{Let }$M^{2}$ \textit{be a helicoidal
surface of first type in }$\mathbb{I}^{3}$ \textit{satisfying }$%
\bigtriangleup r_{i}=\lambda _{i}r_{i},$ $\lambda _{i}\in 
\mathbb{R}
$. \textit{Then it is isotropic minimal and has the form}%
\begin{equation*}
\mathbf{r}\left( u,v\right) =\left( u\cos v,u\sin v,\alpha \ln \left\vert
u\right\vert +hv\right) ,\text{ }\alpha \in 
\mathbb{R}
\backslash \left\{ 0\right\} .
\end{equation*}

\section{Helicoidal surfaces with constant curvature in $\mathbb{I}^{3}$}

Let us consider the helicoidal surface of first type $M^{2}$ in $\mathbb{I}%
^{3}.$ Then the components of the second fundamental form are%
\begin{equation}
\mathfrak{h}_{11}=g^{\prime \prime },\text{ }\mathfrak{h}_{12}=-\frac{h}{u},%
\text{ }\mathfrak{h}_{22}=ug^{\prime }.  \tag{3.1}
\end{equation}%
Thereby, the relative curvature $K$ of $M^{2}$ is%
\begin{equation}
K=\frac{u^{3}g^{\prime }g^{\prime \prime }-h^{2}}{u^{4}}.  \tag{3.2}
\end{equation}%
Assume that $M^{2}$ has constant relative curvature $K_{0}.$ We have two
cases:

\bigskip

\textbf{Case (a).} $K$ vanishes. It follows from $\left( 3.2\right) $ that $%
g^{\prime }g^{\prime \prime }=\frac{h^{2}}{u^{3}}$ or 
\begin{equation}
g^{\prime }\left( u\right) =\left( \alpha -\frac{h^{2}}{u^{2}}\right) ^{%
\frac{1}{2}},\text{ }\alpha \in 
\mathbb{R}
^{+}.  \tag{3.3}
\end{equation}%
After solving $\left( 3.3\right) ,$ we obtain%
\begin{equation*}
g\left( u\right) =\sqrt{\alpha u^{2}-h^{2}}+h\arctan \left( \frac{h}{\sqrt{%
\alpha u^{2}-h^{2}}}\right) .
\end{equation*}

\textbf{Case (b)}. $K$ is a nonzero constant $K_{0}.$ We can rewrite $\left(
3.2\right) $ as%
\begin{equation*}
g^{\prime }g^{\prime \prime }=K_{0}u+\frac{h^{2}}{u^{3}}
\end{equation*}%
or%
\begin{equation}
g^{\prime }\left( u\right) =\left( K_{0}u^{2}-\frac{h^{2}}{u^{2}}+\gamma
\right) ^{\frac{1}{2}},\text{ }\gamma \in 
\mathbb{R}
.  \tag{3.4}
\end{equation}%
By solving $\left( 3.4\right) $, we derive%
\begin{equation*}
g\left( u\right) =\frac{1}{4}\left( 2a\left( u\right) -2h\arctan \left( 
\frac{-2h^{2}+\gamma u^{2}}{2ha\left( u\right) }\right) +\frac{\gamma }{%
\sqrt{K_{0}}}\ln \left\vert \gamma +2\left( K_{0}u^{2}+\sqrt{K_{0}}a\left(
u\right) \right) \right\vert \right) ,
\end{equation*}%
where 
\begin{equation*}
\gamma \in 
\mathbb{R}
,\text{ and }a\left( u\right) =\sqrt{K_{0}u^{4}-h^{2}+\gamma u^{2}}\text{.}
\end{equation*}

Thus, we have the next result

\bigskip

\textbf{Theorem 3.1. }\textit{Let }$M^{2}$ \textit{be a helicoidal surface
in }$\mathbb{I}^{3}$ \textit{with constant relative curvature }$K_{0}.$ 
\textit{Then we have the following items}

\textit{(i) when }$K_{0}=0,$ $M^{2}$ \textit{has the form}%
\begin{equation}
\left\{ 
\begin{array}{l}
\mathbf{r}\left( u,v\right) =\left( u\cos v,u\sin v,g\left( u\right)
+hv\right) , \\ 
g\left( u\right) =\sqrt{\alpha u^{2}-h^{2}}+h\arctan \left( \frac{h}{\sqrt{%
\alpha u^{2}-h^{2}}}\right) ,\text{ }\alpha \in 
\mathbb{R}
^{+},%
\end{array}%
\right.  \tag{3.5}
\end{equation}

(ii) \textit{otherwise, i.e. }$K_{0}\neq 0$\textit{, it is of the form}%
\begin{equation}
\left\{ 
\begin{array}{l}
\mathbf{r}\left( u,v\right) =\left( u\cos v,u\sin v,g\left( u\right)
+hv\right) , \\ 
g\left( u\right) =\frac{1}{4}\left( 2a\left( u\right) -2h\arctan \left( 
\frac{-2h^{2}+\gamma u^{2}}{2ha\left( u\right) }\right) +\frac{\gamma }{%
\sqrt{K_{0}}}\ln \left\vert \gamma +2\left( K_{0}u^{2}+\sqrt{K_{0}}a\left(
u\right) \right) \right\vert \right) , \\ 
a\left( u\right) =\sqrt{K_{0}u^{4}-h^{2}+\gamma u^{2}},\text{ }\gamma \in 
\mathbb{R}
\text{.}%
\end{array}%
\right.  \tag{3.6}
\end{equation}

\bigskip

\textbf{Example 3.2. }Take $h=1$, $\alpha =1$, $u\in \left[ 1,5\right] $ and 
$v\in \left[ 0,4\pi \right] $ in $\left( 3.5\right) .$ Then $M^{2}$ becomes
isotropic flat and can be drawn as in Fig 1.%
\begin{gather*}
\FRAME{itbpF}{1.2151in}{1.8213in}{0in}{}{}{Figure}{\special{language
"Scientific Word";type "GRAPHIC";maintain-aspect-ratio TRUE;display
"USEDEF";valid_file "T";width 1.2151in;height 1.8213in;depth
0in;original-width 1.1813in;original-height 1.7841in;cropleft "0";croptop
"1";cropright "1";cropbottom "0";tempfilename
'O0U1R200.wmf';tempfile-properties "XPR";}} \\
\text{\textbf{Fig. 1. }\textit{A helicoidal surface with }}K_{0}=0,\text{ }%
h=1.
\end{gather*}

\bigskip

The isotropic mean curvature $H$ of $M^{2}$ is given by%
\begin{equation*}
H=\frac{g^{\prime }}{u}+g^{\prime \prime }.
\end{equation*}%
Suppose that $M^{2}$ has constant isotropic mean curvature $H_{0}.$ Then
putting $g^{\prime }=p,$ we obtain the following Riccati equation%
\begin{equation}
p^{\prime }+\frac{p}{u}=H_{0}.  \tag{3.7}
\end{equation}%
Solving $\left( 3.7\right) ,$ we get%
\begin{equation*}
g\left( u\right) =\frac{H_{0}}{4}u^{2}+\alpha \ln \left( u\right) +\beta
\end{equation*}%
for some constants $\alpha ,\beta \in 
\mathbb{R}
$ and $\alpha \neq 0.$

Therefore we have proved the next result:

\bigskip

\textbf{Theorem 3.3. }\textit{Let }$M^{2}$ \textit{be a helicoidal surface
in }$\mathbb{I}^{3}$ \textit{with constant isotropic mean curvature }$H_{0}.$
\textit{Then it has the following form}%
\begin{equation}
\left\{ 
\begin{array}{l}
\mathbf{r}\left( u,v\right) =\left( u\cos v,u\sin v,g\left( u\right)
+hv\right) , \\ 
g\left( u\right) =\frac{H_{0}}{4}u^{2}+\alpha \ln \left( u\right) +\beta ,%
\text{ }\alpha \in 
\mathbb{R}
\backslash \left\{ 0\right\} .%
\end{array}%
\right.  \tag{3.8}
\end{equation}

\bigskip

\textbf{Example 3.4. }Let us put $h=1.5$, $H_{0}=-\alpha =-1,$ $\beta =0,$ $%
u\in \left[ 1,5\right] $ and $v\in \left[ -\pi ,\pi \right] $ in $\left(
3.8\right) .$ Then we draw it as in Fig 2.%
\begin{gather*}
\FRAME{itbpF}{1.9787in}{2.0297in}{0in}{}{}{Figure}{\special{language
"Scientific Word";type "GRAPHIC";maintain-aspect-ratio TRUE;display
"USEDEF";valid_file "T";width 1.9787in;height 2.0297in;depth
0in;original-width 1.9398in;original-height 1.9917in;cropleft "0";croptop
"1";cropright "1";cropbottom "0";tempfilename
'NYD3OO00.wmf';tempfile-properties "XPR";}} \\
\text{\textbf{Fig. 2. }\textit{A helicoidal surface with }}H_{0}=-1,\text{ }%
h=1.5.
\end{gather*}

\section{Special curves on the helicoidal surfaces in $\mathbb{I}^{3}$}

For more details of special curves on the surfaces in $\mathbb{I}^{3}$ such
as, geodesics, asymptotic lines and lines of curvature, see \cite{35}, p.
163-181.

In this section we aim to investigate such curves on a helicoidal surface in 
$\mathbb{I}^{3}.$

Let $M^{2}$ be a helicoidal surface in $\mathbb{I}^{3}$, then any point of a
curve on $M^{2}$ has the position vector%
\begin{equation}
\mathbf{r}\left( u\left( s\right) ,v\left( s\right) \right) =\mathbf{r}%
\left( s\right) =\left( u\left( s\right) \cos \left( v\left( s\right)
\right) ,u\left( s\right) \sin \left( v\left( s\right) \right) ,g\left(
u\left( s\right) \right) +hv\left( s\right) \right) ,  \tag{4.1}
\end{equation}%
where $s$ is arc-length parameter of $\mathbf{r}\left( s\right) $. Denote
the derivative with respect to $s$ by a dot. Then $\mathbf{t}\left( s\right)
=\mathbf{\dot{r}}=\left( t_{1}\left( s\right) ,t_{2}\left( s\right)
,t_{3}\left( s\right) \right) $ is the tangent vector of $\mathbf{r}\left(
s\right) $. We can take a side tangential vector $\sigma =\left( \sigma
_{1}\left( s\right) ,\sigma _{2}\left( s\right) ,\sigma _{3}\left( s\right)
\right) $ in the tangent plane of $M^{2}$ such that 
\begin{equation*}
\sigma _{1}^{2}+\sigma _{2}^{2}=1,\text{ }\sigma _{1}t_{1}+\sigma _{2}t_{2}=0%
\text{, }t_{1}\sigma _{2}-t_{2}\sigma _{1}=1.
\end{equation*}%
Therefore we have an orthonormal triple $\left\{ \mathbf{t},\sigma ,\mathbf{N%
}=\left( 0,0,1\right) \right\} .$ The second derivative of $\mathbf{r}\left(
s\right) $ with respect to $s$ has the following decomposition%
\begin{equation*}
\mathbf{\ddot{r}}=\kappa _{g}\sigma +\kappa _{n}\mathbf{N,}
\end{equation*}%
where $\kappa _{g}$ and $\kappa _{n}$ are respectively called the\textit{\
geodesic curvature} and \textit{normal curvature} of $\mathbf{r}\left(
s\right) $ on $M^{2}.$ The curve $\mathbf{r}\left( s\right) $ is called 
\textit{geodesic }(resp.,\textit{\ asymptotic line}) if and only if its
geodesic curvature $\kappa _{g}$ (resp., normal curvature $\kappa _{n}$)
vanishes.

Also the first derivative of $\sigma \left( s\right) $ with respect to $s$
has the decomposition%
\begin{equation*}
\dot{\sigma}=-\kappa _{g}\mathbf{t+\tau }_{g}\mathbf{N,}
\end{equation*}%
in which $\mathbf{\tau }_{g}$ is called the \textit{geodesic torsion }of $%
\mathbf{r}\left( s\right) $ on $M^{2}.$

In terms of the components of the fundamental forms of $M^{2},$ the side
tangential vector $\sigma $ is given by%
\begin{equation*}
\sigma =-\frac{1}{\det \left( \mathfrak{g}_{ij}\right) }\left[ \left( 
\mathfrak{g}_{12}\dot{u}+\mathfrak{g}_{22}\dot{v}\right) \mathbf{r}%
_{u}-\left( \mathfrak{g}_{11}\dot{u}+\mathfrak{g}_{12}\dot{v}\right) \mathbf{%
r}_{v}\right] .
\end{equation*}%
So, the geodesic curvature of $\mathbf{r}\left( s\right) $ on the helicoidal
surface $M^{2}$ in\textit{\ }$\mathbb{I}^{3}$ is given by%
\begin{equation}
\kappa _{g}\left( s\right) =u^{2}\left( \dot{v}\right) ^{3}-u\dot{u}\ddot{v}%
-2\left( \dot{u}\right) ^{2}\dot{v}-u\dot{v}\ddot{u}.  \tag{4.2}
\end{equation}%
It is easliy from $\left( 4.2\right) $ that the curves $v=const.$ on $M^{2}$
are geodesics of $M^{2}$ but not the curves $u=const.,$ which implies the
next result.

\bigskip

\textbf{Theorem 4.1. }\textit{The }$v-$\textit{parameter curves on the
helicoidal surfaces in }$\mathbb{I}^{3}$\textit{\ are geodesics but not }$u-$%
\textit{parameter curves.}

\bigskip

The normal curvature of $\mathbf{r}\left( s\right) $ on $M^{2}$ in\textit{\ }%
$\mathbb{I}^{3}$ is 
\begin{equation}
\kappa _{n}\left( s\right) =g^{\prime \prime }\left( \dot{u}\right) ^{2}-2%
\frac{h}{u}\left( \dot{u}\dot{v}\right) +ug^{\prime }\left( \dot{v}\right)
^{2}.  \tag{4.3}
\end{equation}%
By $\left( 4.3\right) ,$ the curves $u=const.$ are asymptotic lines of $%
M^{2} $ if and only if $g$ is a constant function. Similarly the curves $%
v=const.$ are asymptotic lines of $M^{2}$ if and only if $g$ is a linear
function.

\bigskip

Hence, we have proved the following

\bigskip

\textbf{Theorem 4.2. }\textit{(i) The }$u-$\textit{parameter curves on a
helicoidal surface in }$\mathbb{I}^{3}$ \textit{are asymptotic curves if and
only if it is a helicoid from Euclidean perspective;}

\textit{(ii) the }$v-$\textit{parameter curves on the helicoidal surfaces in 
}$\mathbb{I}^{3}$\textit{\ are asymptotic curves if and only if }$g$\textit{%
\ is a linear function.}

\bigskip

On the other hand a curve on a surface is called a \textit{line of curvature}
if its geodesic torsion $\tau _{g}$ vanishes. The function $\tau _{g}$ can
be defined as%
\begin{equation*}
\tau _{g}=\frac{%
\begin{vmatrix}
dv^{2} & -dudv & du^{2} \\ 
\mathfrak{g}_{11} & \mathfrak{g}_{12} & \mathfrak{g}_{22} \\ 
\mathfrak{h}_{11} & \mathfrak{h}_{12} & \mathfrak{h}_{22}%
\end{vmatrix}%
}{\det \left( \mathfrak{g}_{ij}\right) \mathfrak{g}}.
\end{equation*}%
Hence,\textbf{\ }a curve on $M^{2}$ in $\mathbb{I}^{3}$ is a line of
curvature if and only if the following equation satisfies%
\begin{equation*}
-\left( \frac{h}{u}\right) \left( \dot{u}\right) ^{2}+\left( ug^{\prime
}-u^{2}g^{\prime \prime }\right) \dot{u}\dot{v}+\left( hu\right) \left( \dot{%
v}\right) ^{2}=0.
\end{equation*}

Therefore we can give the following result.

\bigskip

\textbf{Theorem 4.3. }\textit{The parameter curves on the helicoidal
surfaces in }$\mathbb{I}^{3}$\textit{\ are lines of curvature if and only if
those are surfaces of revolution.}

\end{document}